\documentclass{amsart}
\usepackage{amssymb,latexsym,amsmath,amsfonts}
\usepackage[colorlinks,linkcolor=blue,anchorcolor=blue,citecolor=blue]{hyperref}

\newcommand\C{{\mathbb C}}

\newcommand\Q{{\mathbb Q}}

\newcommand\Z{{\mathbb Z}}

\newcommand\PP{{\mathbb P}}

\newcommand\SL{{\mathrm {SL}}}

\newcommand\HH{{\mathcal H}}
\newcommand\OO{{\mathcal O}}
\newcommand\NN{{\mathcal{N}}}

\newcommand\ns{\mathrm{ns}}

\newcommand\h{\mathrm{h}}

\newtheorem{theorem}{Theorem}[section]

\newtheorem{example}[theorem]{Example}
\numberwithin{equation}{section}

\begin{document}
\title[Bounding the $j$-invariant of integral points]{Bounding the $j$-invariant of integral points on certain modular curves}

\author{Min Sha}
\address{Institut de Mathematiques de Bordeaux, Universite Bordeaux 1
, 33405 Talence Cedex, France}
%\curraddr{Department of Mathematics and Statistics,
%Case Western Reserve University, Cleveland, Ohio 43403}
\email{shamin2010@gmail.com}
\thanks{%The authors were partially supported by the {\it Agence nationale de la recherche} project ``HAMOT".
The author was supported by the China Scholarship Council.}



\date{}

\begin{abstract}
In this paper, we obtain two effective bounds for the $j$-invariant of integral points on certain modular curves which have positive genus and
less than three cusps.
\end{abstract}

\maketitle

\section{Introduction}
Throughout this paper, let $N\ge 2$ be a positive integer and $X(N)$ the principal modular curve of level $N$. Let $\Gamma$ be a congruence subgroup of level $N$ of $\SL_2(\Z)$ and $X_{\Gamma}$ its corresponding modular curve. We denote by $\nu_{\infty}(\Gamma)$ the number of cusps of $X_{\Gamma}$.

Assume that $X_{\Gamma}$ is defined over a number field $K$. Let $S$ be a finite set of absolute values of $K$,
containing all the Archimedean valuations and normalized with respect to $\Q$.
We call a $K$-rational point $P\in X_{\Gamma}(K)$ an $S$-integral point if $j(P)\in \mathcal{O}_{S}$,
where $j$ is the standard $j$-invariant function on $X_{\Gamma}$ and $\mathcal{O}_{S}$ is the ring of $S$-integers in $K$.

By Siegel's theorem~\cite{Si29}, $X_{\Gamma}$ has only finitely many $S$-integral points
when $X_{\Gamma}$ has positive genus or $\nu_{\infty}(\Gamma)\ge 3$. But the existing proofs of Siegel's theorem are not effective, that is they don't
 provide with any effective bounds for the $j$-invariant of $S$-integral points. In 1995, Bilu \cite[Proposition 5.1]{Bi95} showed the following effective result without giving a quantitative version.
 \begin{theorem}[Bilu \cite{Bi95}]\label{Bilu1}
Siegel's theorem is effective for $X_{\Gamma}$ if

1.  $\nu_{\infty}(\Gamma)\ge 3$, or

2. $\Gamma$ has no elliptic elements.
\end{theorem}

Afterwards, combining with Chevalley-Weil theorem, Bilu \cite[Proposition 12]{Bilu02} gave the following refinement. But there were still no quantitative results therein.
\begin{theorem}[Bilu \cite{Bilu02}]\label{Bilu2}
Assume that $\Gamma$ has a congruence subgroup $\Gamma^{\prime}$ with $\nu_{\infty}(\Gamma^{\prime})\ge 3$, and $\Gamma^{\prime}$ contains all elliptic elements of $\Gamma$. Then Siegel's theorem is effective for $X_{\Gamma}$.
\end{theorem}

Most recently, Sha \cite{Sha} obtained a quantitative version for Theorem \ref{Bilu1} when $\nu_{\infty}(\Gamma)\ge 3$. Especially,  Bajolet and Sha \cite{BaSh1} gave much better bounds in the non-split Cartan case. We briefly recall the main results of \cite{Sha} as follows for the convenience of the reader.

Put $d=[K:\Q]$ and $s=|S|$. Let $\OO_K$ be the ring of integers of $K$.
We define the following quantity
\begin{equation}\label{Delta_0}
\Delta_0(N)=d^{-d}\sqrt{N^{dN}|D|^{\varphi(N)}}\left(\log(N^{dN}|D|^{\varphi(N)})\right)^{d\varphi(N)}
\left(\prod\limits_{\substack{v\in S\\v\nmid \infty}}\log\NN_{K/\Q}(v)\right)^{\varphi(N)}
\end{equation}
as a function of $N$, where $D$ is the absolute discriminant of $K$, $\varphi(N)$ is the Euler's totient function, and the norm $\NN_{K/\Q}(v)$ of a place $v$, by definition, is equal to $|\OO_K/\mathfrak{p}_v|$ when $v$ is finite and $\mathfrak{p}_v$ is its corresponding prime ideal, and is set to be 1 if $v$ is infinite.

Let $p$ be the maximal rational prime below $S$, with the convention $p=1$ if $S$ only consists
of the infinite places. We denote by $\h(\cdot)$ the usual absolute logarithmic height.
For $P\in X_{\Gamma}(\bar{K})$, we write $\h(P)=\h(j(P))$, where $\bar{K}$ is a fixed algebraic closure of $K$.

In addition, when $N$ is a prime power, we define $M=3N$ if $N$ is a power of 2, and $M=2N$ otherwise.

\begin{theorem}[Sha \cite{Sha}]\label{main}
Assume that $\nu_{\infty}(\Gamma)\ge 3$.
Then for any $S$-integral point $P$ on $X_{\Gamma}$, the following hold.

{\rm (1)} If $N$ is not a power of any prime, we have
$$
\h(P)\le \left(CdsN^{2}\right)^{2sN}(\log(dN))^{3sN}p^{dN}\Delta_0(N),
$$
where $C$ is an absolute effective constant.

{\rm (2)} If $N$ is a power of some prime, we have
$$
\h(P)\le \left(CdsM^{2}\right)^{2sM}(\log(dM))^{3sM}p^{dM}\Delta_0(M),
$$
where $C$ is an absolute effective constant.
\end{theorem}

In this paper, we will give a quantitative version for Theorem \ref{Bilu2} by applying Theorem \ref{main}. As an application, it is also a quantitative version for Theorem \ref{Bilu1} when $\Gamma$ has no elliptic elements, as well as for certain modular curves which have positive genus and
less than three cusps. For example, the classical modular curve $X_{0}(p)$ for a prime $p>13$ has positive genus and two cusps.

Let $\mathcal{H}$ denote the Poincar$\acute{\rm e}$ upper half-plane: $\mathcal{H}=\{\tau\in\C: {\rm Im} \tau>0\}$.
Recall that a non-cuspidal point $P\in X_\Gamma$ is called \emph{elliptic} if for some $z\in\HH$ representing $P$ the stabilizer $\Gamma_{z}\neq\{\pm 1\}$. Notice that the curve $X_\Gamma$ has finitely many elliptic points. We assume that the set of its elliptic points is
 $\{P_1,P_2,\cdots,P_n\}$. For each elliptic point $P_i$, we fix a pre-image $z_i$ in $\HH$.
 We denote by $\Gamma_{z_i}$ the stabilizer of $z_i$ in $\Gamma$. It is well-known that each $\Gamma_{z_i}$ is cyclic of order 3, 4, or 6.

 Let $\widetilde\Gamma$ be the congruence subgroup generated by the principal congruence subgroup
 $\Gamma(N)$ and $\{\Gamma_{z_1},\cdots,\Gamma_{z_n}\}$. Consider the natural finite covering $\phi: X_{\widetilde\Gamma}\to X_\Gamma$. For any point $\tilde P\in X_{\widetilde\Gamma}$, fix a pre-image $z\in \HH$, the ramification index of $\tilde P$ over $X_\Gamma$ is equal to the index  $[\pm\Gamma_{z}:\pm\widetilde\Gamma_{z}]$ which does not depend on the choice of $z$.
 Therefore, $\phi$ is unramified outside the cusps.

Assume that $\Gamma$ has a congruence subgroup $\Gamma^\prime$ such that $X_{\Gamma^\prime}$ has at least three cusps and the finite covering $X_{\Gamma^\prime}\to X_{\Gamma}$ is unramified outside the cusps. Then we must have $\widetilde\Gamma\subseteq\Gamma^\prime$, subsequently $X_{\widetilde\Gamma}$ also has at least three cusps. Under this assumption, by using Theorem \ref{main} we can get an effective Siegel's theorem for $X_{\widetilde\Gamma}$. Then the effective Siegel's theorem for $X_\Gamma$ follows from quantitative Riemann existence theorem \cite{BS} and quantitative Chevalley-Weil theorem \cite{BSS}.

For stating the main result, we need to fix some other notation. Put
\begin{equation}
d_{N}=\left\{ \begin{array}{ll}
                 \frac{1}{2}N^{3}\prod_{q|N}(1-1/q^{2}) & \textrm{if $N>2$},\\
                  \\
                 6 & \textrm{if $N=2$},
                 \end{array} \right.
\notag
\end{equation}
where $q$ runs through all primes dividing $N$. Define
\begin{equation}
D^{*}=D^{d_N}e^{(\h(S)+(1+\log 1728)\Lambda)dd_N},\notag
\end{equation}
where
$$
\Lambda=\left((\frac{d_{N}(N-6)}{12N}+2)d_N\right)^{25(\frac{d_{N}(N-6)}{12N}+2)d_N},
$$
and
$$
\h(S)=\frac{\sum_{v\in S}\log\NN_{K/\Q}(v)}{d}.
$$
Next we define
\begin{equation}\label{Delta}
\Delta(N)=d^{-d}\sqrt{N^{Ndd_{N}}|D^{*}|^{\varphi(N)}}\left(\log(N^{Ndd_{N}}|D^{*}|^{\varphi(N)})\right)^{\varphi(N)dd_{N}}
\left(\prod\limits_{\substack{v\in S\\v\nmid \infty}}\log\NN_{K/\Q}(v)\right)^{\varphi(N)d_N}.
\end{equation}

  Now we are ready to state the main result.
\begin{theorem}\label{main1}
Assume that $\Gamma$ has a congruence subgroup $\Gamma^{\prime}$ with $\nu_{\infty}(\Gamma^{\prime})\ge 3$, and $\Gamma^{\prime}$ contains all elliptic elements of $\Gamma$.
Then for any $S$-integral point $P$ on $X_{\Gamma}$, the following hold.

{\rm (1)} If $N$ is not a power of any prime, we have
$$
\h(P)\le \left(Cdsd_{N}^{2}N^{2}\right)^{2sNd_N}(\log(dNd_N))^{3sNd_N}p^{dNd_N}\Delta(N),
$$
where $C$ is an absolute effective constant.

{\rm (2)} If $N$ is not a power of any prime, we have
$$
\h(P)\le \left(Cdsd_{M}^{2}M^{2}\right)^{2sMd_M}(\log(dMd_M))^{3sMd_M}p^{dMd_M}\Delta(M),
$$
where $C$ is an absolute effective constant.
\end{theorem}

Here, we would like to give some examples satisfying the assumptions in Theorem \ref{main1}.
\begin{example}
{\rm
Assume that $\Gamma$ has no elliptic elements. Then the principal congruence subgroup $\Gamma(N)$ is such a subgroup of $\Gamma$ when $N\ge 2$.
}
\end{example}

\begin{example}
{\rm
For a prime $p>13$, the classical modular curve $X_{0}(p)$ has positive genus and two cusps.
By \cite[Proof of Theorem 10]{Bilu02}, it has a congruence subgroup $\Gamma^{\prime}$ with $\nu_{\infty}(\Gamma^{\prime})\ge 3$, and $\Gamma^{\prime}$ contains all elliptic elements of $\Gamma_0(p)$.
}
\end{example}

\begin{example}
{\rm
Assume that $\Gamma_{z_1},\cdots,\Gamma_{z_n}$ generate a finite subgroup $G$ and $|G|<\frac{1}{4}N^{2}\prod\limits_{q|N}(1-q^{-2})$,
where the product being taken over all primes $q$ dividing $N$. By \cite[Corollary 2.4]{BI11}, $X_{\widetilde\Gamma}$ has at least three cusps.
Then $\widetilde\Gamma$ is such a subgroup of $\Gamma$.
}
\end{example}

\section{Quantitative Riemann existence theorem for $X_{\widetilde\Gamma}$}
The Riemann Existence Theorem asserts that every
compact Riemann surface is (analytically isomorphic to) a complex algebraic
curve. Bilu and Strambi \cite[Theorem 1.2]{BS} gave a quantitative version of Riemann Existence Theorem, which is a key tool in this paper.

Notice that the $j$-invariant induces naturally two coverings $X_{\Gamma}\to\PP^{1}(\C)$ and
 $X_{\widetilde\Gamma}\to\PP^{1}(\C)$, respectively. We use the same notation $j$ to denote both of them without confusions.
In addition, the $j$-invariant also defines an isomorphism $X(1)\cong \PP^{1}(\C)$.

For the covering $j: X_{\widetilde\Gamma}\to\PP^{1}(\C)$, we assume that its degree is $\tilde n$ and the genus of the curve $X_{\widetilde\Gamma}$ is $\tilde g$. Then by Theorem 1.2 of \cite{BS}, there exists a rational function $y\in \bar{K}(X_{\widetilde\Gamma})$
such that $\bar{K}(X_{\widetilde\Gamma})=\bar{K}(j,y)$ and the rational functions $j,y\in \bar{K}(X_{\widetilde\Gamma})$ satisfy the equation $\tilde f(j,y)=0$, where $\tilde f(X,Y)\in \bar{K}[X,Y]$ is an absolutely irreducible polynomial satisfying
\begin{equation}
\deg_X \tilde f=\tilde g+1,\qquad \deg_Y \tilde f=\tilde n.
\end{equation}

Consider the natural sequence of coverings $X(N)\to X_{\widetilde\Gamma}\to\PP^{1}(\C)$. Applying the formula in the bottom of \cite[Page 101]{Diamond}, we know that the degree of the covering $X(N)\to\PP^{1}(\C)$ is $d_N$.
Combining with the genus formula of $X(N)$ (see \cite[Page 108, Figure 3.4]{Diamond}), we have
\begin{equation}\label{bound}
\tilde n\le d_{N},\qquad \tilde g\le 1+\frac{d_{N}(N-6)}{12N}.
\end{equation}

\section{Quantitative Chevalley-Weil theorem for $\phi:X_{\widetilde\Gamma}\to X_\Gamma$}
The Chevalley-Weil theorem asserts that for an \'etale covering of projective varieties over a number field $F$, the discriminant of the field of definition of the fiber over an $F$-rational point is uniformly bounded. Bilu, Strambi and Surroca \cite{BSS} got a fully explicit version of this theorem in dimension one, which is another key tool of this paper.

The covering $j: X_{\Gamma}\to\PP^{1}(\C)$ is unramified above all affine points except 0 and 1728 which correspond to elliptic elements of $\SL_2(Z)$.
For the covering $\phi: X_{\widetilde\Gamma}\to X_\Gamma$, it is unramified outside the cusps. Notice that the poles of $j$ are exactly the cusps. Then by Theorem 1.6 of \cite{BSS}, for every
$P\in X_{\Gamma}(K)$ and $\tilde P\in X_{\widetilde\Gamma}(\bar{K})$ such that $\phi(\tilde P)=P$, we have
\begin{equation}
\mathcal{N}_{K/\Q}(D_{K(\tilde P)/K})\le e^{[K(\tilde P):\Q]\cdot(\h(S)+(1+\log 1728)\tilde\Lambda)},
\end{equation}
where $D_{K(\tilde P)/K}$ is the relative discriminant of $K(\tilde P)/K$, and $\tilde\Lambda=((\tilde g+1)\tilde n)^{25(\tilde g+1)\tilde n}$. According to (\ref{bound}), we have $\tilde\Lambda\le \Lambda$.
Hence
\begin{equation}\label{disc}
\mathcal{N}_{K/\Q}(D_{K(\tilde P)/K})\le e^{[K(\tilde P):\Q]\cdot(\h(S)+(1+\log 1728)\Lambda)}.
\end{equation}

Notice that $[K(\tilde P):K]=[K(\tilde P):K(P)]$ is not greater than the degree of $\phi$.
So by (\ref{bound}), we have
\begin{equation}
[K(\tilde P):K]\le d_N.
\end{equation}

\section{Proof of Theorems}
Under the assumptions of Theorem \ref{main1}, the curve $X_{\widetilde\Gamma}$ has at least three cusps. In this section, we fix an $S$-integral point $P$ on $X_\Gamma$ and a point $\tilde P$ on $X_{\widetilde\Gamma}$ such that $\phi(\tilde P)=P$.

Let $K_0=K(\tilde P)$ and $d_0=[K_0:\Q]$. Let $S_0$ be the set consisting of the extensions of the places from $S$ to $K_{0}$, that is,
$$
S_0=\{v\in M_{K_0}:v|w\in S\},
$$
where $M_{K_{0}}$ is the set of all valuations (or places) of $K_{0}$ extending the standard infinite and $p$-adic valuations of $\Q$.
Put $s_0=|S_0|$. We define the following quantity
\begin{equation}
\Delta_1=d_{0}^{-d_{0}}\sqrt{N^{d_{0}N}|D_{0}|^{\varphi(N)}}\left(\log(N^{d_{0}N}|D_{0}|^{\varphi(N)})\right)^{d_{0}\varphi(N)}
\left(\prod\limits_{\substack{v\in S_{0}\\v\nmid \infty}}\log\NN_{K_{0}/\Q}(v)\right)^{\varphi(N)},
\notag
\end{equation}
where $D_{0}$ is the absolute discriminant of $K_{0}$.

Notice that $d_0\le dd_N$ and $s_0\le sd_{N}$. Let $D_{K_{0}/K}$ be the relative discriminant of $K_{0}/K$.
By (\ref{disc}), we have
\begin{align*}
D_0&=\mathcal{N}_{K/\Q}(D_{K_{0}/K})D^{[K_{0}:K]}\\
&\le D^{*}.
\end{align*}

Now let $w$ be a non-archimedean place of $K$, and $v_{1},\cdots,v_{m}$ all its
extensions to $K_0$, their residue degrees over $K$ being $f_{1},\cdots,f_{m}$ respectively.
Then $f_{1}+\cdots+f_{m}\le [K_{0}:K]\le d_N$, which implies that $f_{1}\cdots f_{m}\le 2^{d_N}$.
Since $\NN_{K_0/\Q}(v_{k})=\NN_{K/\Q}(w)^{f_{k}}$ for $1\le k\le m$, we have
$$
\prod\limits_{v|w}\log\NN_{K_{0}/\Q}(v)\le 2^{d_N}(\log\NN_{K/\Q}(w))^{d_N}.
$$
Hence
\begin{equation}\label{logN}
\prod\limits_{\substack{v\in S_{0}\\v\nmid \infty}}\log\NN_{K_{0}/\Q}(v)
\le 2^{sd_N}\left(\prod\limits_{\substack{v\in S\\v\nmid \infty}}\log\NN_{K/\Q}(v)\right)^{d_N}.
\notag
\end{equation}

Combining with $d_0\ge d$, we have
$$\Delta_1\le 2^{s\varphi(N)d_N}\Delta(N).$$

First we assume that $N$ is not a power of any prime. By Theorem \ref{main} (1), we have
$$
\h(\tilde P)\le \left(Cd_{0}s_{0}N^{2}\right)^{2s_{0}N}(\log(d_{0}N))^{3s_{0}N}p^{d_{0}N}\Delta_1,
$$
where $C$ is an absolute effective constant. Note that $j(P)=j(\tilde P)$, we have $\h(P)=\h(\tilde P)$. Then we have
\begin{equation}
\h(P)\le \left(Cdsd_{N}^{2}N^{2}\right)^{2sNd_N}(\log(dNd_N))^{3sNd_N}p^{dNd_N}\Delta(N),
\end{equation}
the constant $C$ being modified. So we prove Theorem \ref{main1} (1).

For the case that $N$ is a prime power, applying Theorem \ref{main} (2), we can easily prove
Theorem \ref{main1} (2).

\section*{Acknowledgement}
The author would like to thank his advisor Yuri Bilu for helpful discussions and valuable suggestions. He is also grateful to the referee for careful reading and very useful comments.


\begin{thebibliography}{99}

\bibitem{BaSh1}
A. Bajolet and M. Sha, Bounding the $j$-invariant of integral points  on $X_{\ns}^{+}(p)$, \emph{Proceedings of the American Mathematical Society}, to appear. arXiv:1203.1187

\bibitem{Bilu02}
Yu. Bilu, Baker's method and modular curves, \emph{A Panorama of Number Theory or The View from Baker's Garden}
(edited by G. W\"ustholz), 73-88, Cambridge University Press, 2002.

\bibitem{BS}
Yu. Bilu and M. Strambi, Quantitative Riemann existence theorem over a number field, \emph{Acta Arithmetica} \textbf{145} (2010), 319-339.

\bibitem{BSS}
Yu. Bilu, M. Strambi and A. Surroca, Quantitative Chevalley-Weil Theorem for Curves, \emph{Monatshefte f\"ur Mathematik}, to appear. arXiv:0908.1233

\bibitem{Diamond}
F. Diamond and J. Shurman, \emph{A First Course in Modular Forms}, Springer, New York, 2005.


\bibitem{Sha}
M. Sha, Bounding the $j$-invariant of integral points on modular curves, \emph{International Mathematics Research Notices}, doi: 10.1093/imrn/rnt085. arXiv:1208.1337


\bibitem{Bi95}
Yu. Bilu, Effective analysis of integral points on algebraic curves, \emph{Israel Joural of Mathematics} \textbf{90} (1995), 235-252.

\bibitem{BI11}
Yu. Bilu and M. Illengo, Effective Siegel's Theorem for Modular Curves, \emph{Bulletin of the London Mathematical Society}  \textbf{43} (2011), 673-688.

\bibitem{Si29}
C.L. Siegel, \emph{$\ddot{\rm U}$ber einige Anwendungen diophantischer Approximationen},\emph{Abhandlungen der Preussischen Akademie der Wissenschaften}  (1929), Nr. 1. (=Ges. Abh. I, 209-266, Springer, 1966.)

\end{thebibliography}
\end{document}